# Easy and nearly simultaneous proofs of the Ergodic Theorem and Maximal Ergodic Theorem

**Michael Keane[1] and Karl Petersen[2]**

*Wesleyan University and University of North Carolina*

**Abstract:** We give a short proof of a strengthening of the Maximal Ergodic Theorem which also immediately yields the Pointwise Ergodic Theorem.

Let $(X, \mathcal{B}, \mu)$ be a probability space, $T : X \to X$ a (possibly noninvertible) measure-preserving transformation, and $f \in L^1(X, \mathcal{B}, \mu)$. Let

$$A_k f = \frac{1}{k} \sum_{j=0}^{k-1} fT^j, \quad f_N^* = \sup_{1 \le k \le N} A_k f, \quad f^* = \sup_N f_N^*, \quad \text{and } \overline{A} = \limsup_{k \to \infty} A_k f.$$

When $\lambda$ is a constant, the following result is the Maximal Ergodic Theorem. Choosing $\lambda = \overline{A} - \epsilon$ covers most of the proof of the Ergodic Theorem.

**Theorem.** *Let $\lambda$ be an invariant ($\lambda \circ T = \lambda$ a.e.) function on $X$ with $\lambda^+ \in L^1$. Then*

$$\int_{\{f^* > \lambda\}} (f - \lambda) \ge 0.$$

*Proof.* We may assume that $\lambda \in L^1\{f^* > \lambda\}$, since otherwise

$$\int_{\{f^* > \lambda\}} (f - \lambda) = \infty \ge 0.$$

But then actually $\lambda \in L^1(X)$, since on $\{f^* \le \lambda\}$ we have $f \le \lambda$, so that on this set $\lambda^- \le -f + \lambda^+$, which is integrable.

Assume first that $f \in L^\infty$. Fix $N = 1, 2, \ldots,$ and let

$$E_N = \{f_N^* > \lambda\}.$$

Notice that

$$(f - \lambda)\chi_{E_N} \ge (f - \lambda),$$

since $x \notin E_N$ implies $(f - \lambda)(x) \le 0$. Thus for a very large $m \gg N$, we can break up

$$\sum_{k=0}^{m-1} (f - \lambda)\chi_{E_N}(T^k x)$$

into convenient strings of terms as follows. There is maybe an initial string of 0's during which $T^k x \notin E_N$. Then there is a first time $k$ when $T^k x \in E_N$, which initiates a string of no more than $N$ terms, the sum of which is positive (using on each of these terms the fact that $(f - \lambda)\chi_{E_N} \geq (f - \lambda)$). Beginning after the last term in this string, we repeat the previous analysis, finding maybe some 0's until again some $T^k x \in E_N$ initiates another string of no more than $N$ terms and with positive sum. The full sum of $m$ terms may end in the middle of either of these two kinds of strings (0's, or having positive sum). Thus we can find $j = m - N + 1, \dots, m$ such that

$$\sum_{k=0}^{m-1} (f - \lambda)\chi_{E_N}(T^k x) \geq \sum_{k=j}^{m-1} (f - \lambda)\chi_{E_N}(T^k x) \geq -N(\|f\|_\infty + \lambda^+(x)).$$

Integrating both sides, dividing by $m$, and letting $m \to \infty$ gives

$$m \int_{E_N} (f - \lambda) \geq -N(\|f\|_\infty + \|\lambda^+\|_1),$$

$$\int_{E_N} (f - \lambda) \geq \frac{-N}{m}(\|f\|_\infty + \|\lambda^+\|_1),$$

$$\int_{E_N} (f - \lambda) \geq 0.$$

Letting $N \to \infty$ and using the Dominated Convergence Theorem concludes the proof for the case $f \in L^\infty$.

To extend to the case $f \in L^1$, for $s = 1, 2, \dots$ let $\phi_s = f \cdot \chi_{\{|f| \leq s\}}$, so that $\phi_s \in L^\infty$ and $\phi_s \to f$ a.e. and in $L^1$. Then for fixed $N$

$$(\phi_s)_N^* \to f_N^* \quad \text{a.e. and in } L^1 \quad \text{and} \quad \mu\left(\{(\phi_s)_N^* > \lambda\} \triangle \{f_N^* > \lambda\}\right) \to 0.$$

Therefore

$$0 \leq \int_{\{(\phi_s)_N^* > \lambda\}} (\phi_s - \lambda) \to \int_{\{f_N^* > \lambda\}} (f - \lambda),$$

again by the Dominated Convergence Theorem. The full result follows by letting $N \to \infty$. $\qquad \square$

**Corollary (Ergodic Theorem).** *The sequence* $(A_k f)$ *converges a.e.*

*Proof.* It is enough to show that

$$\int \overline{A} \leq \int f.$$

For then, letting $\underline{A} = \liminf A_k f$, applying this to $-f$ gives

$$-\int \underline{A} \leq -\int f,$$

so that

$$\int \overline{A} \leq \int f \leq \int \underline{A} \leq \int \overline{A},$$

and hence

$$\int (\overline{A} - \underline{A}) = 0, \quad \text{so that } \overline{A} = \underline{A} \text{ a.e.}$$

Consider first $f^+$ and its associated $\overline{A}$, denoted by $\overline{A}(f^+)$. For any invariant function $\lambda < \overline{A}(f^+)$ such that $\lambda^+ \in L^1$, for example $\lambda = \overline{A}(f^+) \wedge n - 1/n$, we have $\{(f^+)^* > \lambda\} = X$, so the Theorem gives

$$\int f^+ \geq \int \lambda \nearrow \int \overline{A}(f^+).$$

Thus $(\overline{A})^+ \leq \overline{A}(f^+)$ is integrable (and, by a similar argument, so is $(\overline{A})^- \leq \overline{A}(f^-)$.)

Now let $\epsilon > 0$ be arbitrary and apply the Theorem to $\lambda = \overline{A} - \epsilon$ to conclude that

$$\int f \geq \int \lambda \nearrow \int \overline{A}. \qquad \square$$

**Remark.** This proof may be regarded as a further development of one given in a paper by Keane [10], which has been extended to deal also with the Hopf Ratio Ergodic Theorem [8] and with the case of higher-dimensional actions [11], and which was itself a development of the Katznelson-Weiss proof [9] based on Kamae's nonstandard-analysis proof [7]. (It is presented also in the Bedford-Keane-Series collection [1].) Our proof yields both the Pointwise and Maximal Ergodic Theorems essentially simultaneously without adding any real complications. Roughly contemporaneously with this formulation, Roland Zweimüller prepared some preprints [21, 22] also giving short proofs based on the Kamae-Katznelson-Weiss approach, and recently he has also produced a simple proof of the Hopf theorem [23]. Without going too deep into the complicated history of the Ergodic Theorem and Maximal Ergodic Theorem, it is interesting to note some recurrences as the use of maximal theorems arose and waned repeatedly. After the original proofs by von Neumann [18], Birkhoff [2], and Khinchine [12], the role and importance of the Maximal Lemma and Maximal Theorem were brought out by Wiener [19] and Yosida-Kakutani [20], making possible the exploration of connections with harmonic functions and martingales. Proofs by upcrossings followed an analogous pattern. It also became of interest, for instance to allow extension to new areas or new kinds of averages, again to prove the Ergodic Theorem without resort to maximal lemmas or theorems, as in the proof by Shields [16] inspired by the Ornstein-Weiss proof of the Shannon-McMillan-Breiman Theorem for actions of amenable groups [14], or in Bourgain's proofs by means of variational inequalities [3]. Sometimes it was pointed out, for example in the note by R. Jones [6], that these approaches could also with very slight modification prove the Maximal Ergodic Theorem. Of course there are the theorems of Stein [17] and Sawyer [15] that make the connection explicit, just as the transference techniques of Wiener [19] and Calderón [4] connect ergodic theorems with their analogues in analysis like the Hardy-Littlewood Maximal Lemma [5]. In many of the improvements over the years, ideas and tricks already in the papers of Birkhoff, Kolmogorov [13], Wiener, and Yosida-Kakutani have continued to play an essential role.

## Acknowledgment

This note arose out of a conversation between the authors in 1997 at the Erwin Schrödinger International Institute for Mathematical Physics in Vienna, and we thank that institution for its hospitality. Thanks also to E. Lesigne and X. Méla for inducing several clarifications.